\documentclass[10pt]{amsart}
\usepackage{amscd}
\usepackage{amssymb}
\usepackage{latexsym}
\usepackage{graphicx}
\usepackage[parfill]{parskip} 
\usepackage{amsmath}
\usepackage{amsfonts}
\usepackage{amssymb}
\usepackage{latexsym}
\usepackage{graphicx}
\usepackage{enumerate}
\usepackage{bbm}
\usepackage[final]{hyperref}
 \hypersetup{colorlinks=true, linkcolor=blue,
anchorcolor=blue, citecolor=red, filecolor=blue, menucolor=blue,
pagecolor=blue, urlcolor=blue}

\newcommand{\abs}[1]{\vert #1 \vert}

\newcommand{\Bigabs}[1]{\Bigl\vert #1 \Bigr\vert}

\newcommand{\norm}[1]{\left\Vert #1 \right\Vert}

\newcommand{\bignorm}[1]{\bigl\Vert #1 \bigr\Vert}

\newcommand{\R}{\mathbb{R}}

\newcommand{\angles}[1]{\langle #1 \rangle}

 \DeclareMathOperator{\diag}{diag}

\newtheorem{theorem}{Theorem}[section]

\newtheorem{lemma}[theorem]{Lemma}

\theoremstyle{definition}

\theoremstyle{remark}

\numberwithin{equation}{section} 

\title[Almost critical LWP for Monopole equation]{ Almost critical local well-posedness for the space-time Monopole equation in Lorenz gauge}

\author[Achenef Tesfahun]{Achenef Tesfahun}
\address{Department of Mathematical Sciences\\ Norwegian University of Science and Technology\\ Alfred Getz' vei 1\\ N-7491 Trondheim\\ Norway}
\email{tesfahun@math.ntnu.no}

\begin{document}

\maketitle

\begin{abstract}
Recently,  Candy and Bournaveas proved local well-posedness of the space-time monopole equation in Lorenz gauge for initial data in $H^s $ with $s>\frac14$. The equation is $L^2$-critical, and hence a $\frac14$ derivative gap is left between their result and the scaling prediction. In this paper, we consider initial data in the Fourier-Lebesgue space  $\widehat{H_p^s}$ for $1<p\le 2$ which coincides with $H^s$ when $p=2$ but scales like lower regularity Sobolev spaces for $1<p< 2$. In particular, we will see that as $p\rightarrow 1^+$, the critical exponent $s^c_p\rightarrow 1^-$, in which case $\widehat{\dot H_{1+}^{1-}}$ is the critical space. We shall prove almost optimal local well-posedness to the space-time monopole equation in Lorenz gauge with initial data in the aforementioned spaces that correspond to $p$ close to 1.

\end{abstract}

\section{Introduction}\label{section_1}
The space-time Monopole equation can be derived by a dimensional reduction from the Anti-Self-Dual Yang Mills equations, and is given by 
\begin{equation}
\label{MonopEq}
F_A=*D_A\phi,
\end{equation}
where $F_A$ is the curvature of a one-form connection 
$A=A_\alpha dx^\alpha$  ($\alpha=0,1,2$),  $D_A\phi$ is a covariant derivative of the Higgs field  $\phi$ and $*$
is the Hodge star operator with respect to the Minkowski metric $\diag(-1,1,1)$ on $\R^{1+2}$. 
The space-time Monopole  equation was first introduced by Ward ~\cite{w1989} as  a space-time analog of Bogomolny Equations or Magnetic Monopole equations. For more detailed survey of the equation  see ~\cite{dtu2006}.

The unknowns are 
$A_\alpha$ and $\phi$ which are maps from $\R^{1+2}$ into a Lie algebra $\mathfrak{g}$  of a Lie group $G$, which for simplicity we assume to be the matrix group $SU(n)$ :
$$
A_\alpha, \phi : \R^{1+2}\rightarrow \mathfrak{g}.
$$
The curvature of the connection,  $F_A$,  and the covariant derivative of the Higgs field,   $D_A\phi$,  are given by
\begin{align*}
    F_A&=\frac12\left(\partial_\alpha A_\beta-\partial_\beta A_\alpha + [A_\alpha,  A_\beta]\right) d x^\alpha \wedge d x^\beta,
    \\
    D_A\phi&=(\partial_\alpha\phi + [A_\alpha, \phi])dx^\alpha,
    \end{align*}
   where $[\cdot, \cdot]$ is a Lie bracket. Using the definition of the Hodge star operator we write
   $$
   *D_A\phi=-(\partial_t\phi + [A_0, \phi])dx^1 \wedge d x^2-(\partial_1\phi + [A_1, \phi])dx^0\wedge d x^2+ (\partial_2\phi + [A_2, \phi])dx^0\wedge d x^1,
   $$
  which can be equated with $F_A$ to rewrite \eqref{MonopEq} as
   \begin{equation}
   \label{MonopEq2}
   \begin{split}
  \partial_t\phi+ \partial_1 A_2-\partial_2 A_1&= [A_2, A_1]+[\phi, A_0]\\
  \partial_t A_1-\partial_1A_0-\partial_2 \phi&= [A_1, A_0]+[A_2, \phi]\\
  \partial_t A_2 -\partial_2 A_0+\partial_1 \phi&= [A_2, A_0]+[\phi, A_1].
   \end{split}
   \end{equation}
These equations are invariant under the gauge transformations
\begin{equation*}\label{GaugeTransformation}
  A_\alpha\to A_\alpha' =OAO^{-1}+ O\partial_\alpha O^{-1},
  \qquad
  \phi \to \phi' =O\phi O^{-1},
  \end{equation*}
  where 
  $O: \R^{1+2}\rightarrow G$ is a smooth and compactly supported map. The most popular gauges are $(i)$ Temporal gauge: $A_0=0$,  $(ii)$ Coulomb gauge: $\partial^jA_j=0$ and $(iii)$ Lorenz gauge: $\partial^\alpha A_\alpha=0$.
  
 In  Coulomb gauge, the system \eqref{MonopEq2} can be written as  a system of nonlinear wave equations for $(A_1, A_2, \phi)$ that contains null structure in the bilinear terms, coupled with elliptic equation for $A_0$ .  Czubak ~\cite{c2010} proved local well-posedness of the space-time monopole equation \eqref{MonopEq2} in  the Coulumb gauge for small initial data (for both $A$ and $\phi$) in $H^s$ with $s>\frac14$. 
Recently, Candy and Bournaveas ~\cite{cb2012} also observed that the space-time monopole equation in Lorenz gauge can be written as a system of nonlinear wave equations for $(A,  \phi)$ and that   the bilinear terms contained are null forms.  They combined this new structure with bilinear estimates for the homogeneous wave equations in ~\cite{fk2000} to show that the Cauchy problem for \eqref{MonopEq2} in Lorenz gauge  is locally well-posed for large initial data in $H^s$ with $s>\frac14$ . This lifts the smallness assumption on the initial data by Czubak.

On the other hand, the space-time monopole equation is invariant under the scaling 
\begin{equation}
\label{Scaling}
  A_\alpha^\lambda(t,x)=\lambda A_\alpha(\lambda t, \lambda x),
  \qquad
\phi^\lambda(t,x)=\lambda \phi(\lambda t, \lambda x) .
 \end{equation}
Then we have (similarly for $A_\alpha$, since both of the fields scale the same)
$$ \norm{\phi^\lambda(0, x)}_{\dot H^s}=\lambda^{s}
\norm{\phi(0, x)}_{\dot H^s}, $$
which suggests that \eqref{MonopEq2} is $L^2$--critical. So there is still  a $\frac14$ derivative gap between the scaling prediction and the result by Candy--Bournaveas.

In this paper,  instead of the standard $H^s$ space we shall consider initial data in the more general space $\widehat{H_p^s}$ , defined by the norm 
$$
\norm{f}_{\widehat{H_p^s} }= \bignorm{\angles{\xi}^s\hat f(\xi) }_{ L^{p'}_{\xi} }, \quad \frac1p+\frac1{p'}=1, \quad 1<p\le 2, 
$$
where 
$\angles{\cdot}=\sqrt{1+\abs{\cdot}^2}$. The homogeneous version 
$\widehat{\dot H_p^s}$ can be defined similarly. In the special case $p=2$, we have $\widehat{H_p^s}=H^s$. These spaces have been used   by many authors to improve regularity results for a number of dispersive and wave equations. For instance, Gr\"{u}nrock  ~\cite{gr2011} used these spaces to prove almost critical local well-posedness for 3D wave equations with quadratic nonlinearities.  More recently, Grigoryan and Nahmod ~\cite{ga2013} used them to prove almost critical local well-posedness for 2D wave equations with quadratic null forms of Klainerman type.

For the scaling \eqref{Scaling}, it is easy to see  
$$
\norm{\phi^\lambda(0, x)}_{\widehat{\dot H_p^s}}=\lambda^{s+1-\frac2p}
\norm{\phi(0, x)}_{\widehat{\dot H_p^s}},
$$
which suggests that  \eqref{MonopEq2} is $\widehat{\dot H_p^{s^c_p}}$--critical,  where 
$$s^c_p=\frac2p -1.$$
Again, for $p=2$ we have the $L^2$-criticality but as $p\rightarrow 1^+$, $s^c_p\rightarrow 1^-$. 

Our goal in this paper is to prove local well-poesdness of the space-time monopole equation \eqref{MonopEq2} in Lorenz gauge for initial data in  $\widehat{ H_p^{s}}$ for $1<p\le 2$ and $s\ge \frac1p$. Note that $0\le s^c_p< \frac1p$ for  $1<p\le 2$ but $ s^c_p \rightarrow \frac1p$ as $ p\rightarrow 1^+$.  It means that on one end (when $p=1+$) we have almost critical well-posedness in $\widehat{ H_{1+}^{1-}}$ and on the other end (when $p=2$) we have well-posedness in $H^\frac12$ (which is weaker than the result by Candy-Bournaveas).

 To this end,  we complement \eqref{MonopEq2} with initial data
    \begin{equation}
\label{Data}    
     A_\alpha(0)=a_\alpha \in \widehat{H_p^s} , \quad \phi(0)=\phi_0\in \widehat{H_p^s},
   \end{equation}
and state our main result as follows.
  \begin{theorem}
 \label{MainThm}
 Assume $1<p\le 2$ and $s\ge \frac1p$ . Then given initial data \eqref{Data}, there exists $T>0$ and a solution $(A_\alpha, \phi)$ to the space-time monopole equation \eqref{MonopEq2} in Lorenz gauge with regularity 
$$(A_\alpha, \phi)\in C\left([-T, T], \widehat{ H_p^{s}} \right).$$ 
 Moreover, the solution is unique in a certain subspace of this regularity class,  the solution depends continuously on the data, and higher regularity data persists in time. 
  \end{theorem}

The rest of the paper is organized as follows. In Section 2 we reformulate Theorem
\ref{MainThm} to Theorem \ref{Thm} after rewriting \eqref{MonopEq2} in Lorenz gauge. In the same Section we introduce the spaces we shall work in together with some of their properties, and state a Lemma about an estimate in these spaces to the solution of the inhomogeneous linear wave equation.  
In Section 3 we show that  Theorem \ref{Thm} reduces to proving null form estimates. In Section 4 and 5 we prove these null form estimates.  Finally, in Section 6 we give the proof of the Lemma mentioned.

 \section{Rewrititing \eqref{MonopEq2} in Lorenz gauge and restating Theorem \ref{MainThm} }
 We follow ~\cite{cb2012} to rewrite \eqref{MonopEq2} in Lorenz gauge.  Define 
  $u, v: \R^{1+2}\rightarrow \mathfrak{g}\times \mathfrak{g}$ by 
 $$
u=\left(\begin{matrix}
u_1\\
u_2
\end{matrix}\right )=\left(\begin{matrix}
A_0+A_1\\
\phi+A_2
\end{matrix}\right ), \quad v=\left(\begin{matrix}
v_1\\
v_2
\end{matrix}\right )=\left(\begin{matrix}
A_0-A_1\\
\phi-A_2
\end{matrix} 
\right ).
 $$
 so that
 $$
 (A_0, A_1, A_2, \phi)=\frac12 (u_1+ v_1, u_1-v_1, u_2-v_2, u_2+v_2).
  $$
  The Lorenz gauge, $\partial_tA_0-\partial_1A_1-\partial_2 A_2=0, $ becomes 
  $$
  \partial_t(u_1+ v_1)-\partial_1( u_1-v_1)-\partial_2( u_2-v_2)=0.
  $$
  Observe that
     \begin{align*}
 &[A_2, A_1]+[\phi, A_0]\pm ( [A_2, A_0]+[\phi, A_1])=[\phi\pm A_2,  A_0\pm A_1],\\
&[A_1, A_0]+ [A_2, \phi]=\frac12\left( [A_1-A_0, A_1+ A_0]  + [A_2-\phi, A_2+ \phi] \right).
   \end{align*}
Using these identities and the Lorenz gauge, we can rewrite \eqref{MonopEq2} as
  \begin{equation*}
   \begin{split}
  \partial_t u_1- \partial_1 u_1-\partial_2 u_2&= \frac12(u\cdot v-v\cdot u)
  \\
    \partial_t u_2+ \partial_1 u_2-\partial_2 u_1&= [u_2, u_1]
    \\
  \partial_t v_1+\partial_1 v_1+\partial_2 v_2&= \frac12(v\cdot u- u\cdot v)
\\
  \partial_t v_2- \partial_1 v_2+\partial_2 v_1&= [v_2,v_1].
   \end{split}
   \end{equation*}
   Furthermore, by introducing the matrices 
   $$
\alpha_1=\left(\begin{matrix}
1& 0\\
0 &-1
\end{matrix}\right ), \quad \alpha_2=\left(\begin{matrix}
0& 1\\
1& 0
\end{matrix} 
\right ), \quad 
\beta=\left(\begin{matrix}
0& 1\\
-1 &0
\end{matrix}\right ),
 $$
 we can rewrite the above equations as 
  \begin{equation}
   \label{MonopEq2-v2}
   \left\{
\begin{aligned}
  \partial_t u-\alpha\cdot \nabla u&= N(u,v)
  \\
    \partial_t v+\alpha\cdot \nabla v&= N(v,u),
   \end{aligned}
\right.
   \end{equation}
 where $\alpha=(\alpha_1, \alpha_2)$ and 
 $$
 N(u,v)=\left(\begin{matrix}
 \frac12(u\cdot v-v\cdot u)\\
\beta u\cdot u
\end{matrix}\right ).
 $$
 The initial data for \eqref{MonopEq2-v2} become
\begin{equation}
\label{Data2}
u(0)=\left(\begin{matrix}
a_0+a_1\\
\phi_0+a_2
\end{matrix}\right ) , \quad v(0)=\left(\begin{matrix}
a_0-a_1\\
\phi_0-a_2
\end{matrix} 
\right ).
 \end{equation}
 
  We diagonalize \eqref{MonopEq2-v2} by defining the
projections 
$$ \mathcal{P}_{\pm}(\xi) = \frac{1}{2} \left( I \pm \alpha\cdot  \hat\xi  \right), $$
where $ \hat\xi \equiv \frac{\xi}{\abs{\xi}}$. Then $u$ and $v$ split into $u = u_+  + u_-$,  $v = v_+  + v_-$, where $u_\pm = \mathcal{P}_\pm(D) u$ and  $v_\pm = \mathcal{P}_\pm(D) v$. Here $D = -i\nabla$ with symbol $\xi$. Now, applying
$\mathcal{P}_\pm(D)$ to \eqref{MonopEq2-v2}, and using the
identities $\alpha \cdot D = \abs{D} \mathcal{P}_+(D)  - \abs{D}\mathcal{P}_-(D) $,
$\mathcal{P}^2_\pm(D) = \mathcal{P}_\pm(D) $ and $ \mathcal{P}_\pm(D) \mathcal{P}_\mp(D)=0$, 
\eqref{MonopEq2-v2} becomes
  \begin{equation}\label{MonopEq2-v3}
\left\{
\begin{aligned}
  & \bigl(i\partial_t \pm \abs{D} \bigr) u_\pm =i \mathcal{P}_\pm(D) N(u,v),
    \\
  & \bigl( i\partial_t \pm \abs{D}\bigr) v_\pm  =i\mathcal{P}_\pm(D) N(v,u).
\end{aligned}
\right.
\end{equation}
 The initial data for \eqref{MonopEq2-v3} become
\begin{equation}
\label{Data2}
u_\pm(0)=f_\pm \in \widehat{H_p^s} , \quad v_\pm(0)=g_\pm \in \widehat{H_p^s},
 \end{equation}
 where
 \begin{equation*}
f_\pm =\mathcal{P}_\pm(D) \left(\begin{matrix}
a_0+a_1\\
\phi_0+a_2
\end{matrix}\right ), \quad g_\pm=\mathcal{P}_\pm(D)\left(\begin{matrix}
a_0-a_1\\
\phi_0-a_2
\end{matrix}
\right ).
 \end{equation*}
 
 We now reformulate Theorem \ref{MainThm} as follows.
  \begin{theorem}
 \label{Thm}
 Assume $1<p\le 2$ and $s\ge \frac1p$ . Then given initial data \eqref{Data2}, there exists $T>0$ and a solution $(u_\pm, v_\pm)$  to \eqref{MonopEq2-v3} with regularity
$$(u_\pm, v_\pm)\in C\left([-T, T], \widehat{ H_p^{s}} \right).$$   Moreover, the solution is unique in a certain subspace of this regularity class,  the solution depends continuously on the data, and higher regularity data persists in time. 
  \end{theorem}
 
The rest of the paper is dedicated to the proof of Theorem \ref{Thm}. 

Let us first fix some notation and introduce the spaces we shall work in together with some of their properties. 

In estimates we use $A \lesssim B$ as shorthand for $A \le CB$,
where $C\gg 1$ is a positive constant. We use the shorthand
$A\approx B$ for $A\lesssim B\lesssim A$ . Throughout the paper, we use $p'$ as a conjugate exponent for $p$, and $\mathbbm{1}_{\{ \cdot\}}$ denotes the indicator function which is 1 if the condition in the bracket is satisfied and 0 otherwise.

Consider now the inhomogeneous linear equation
\begin{equation}\label{LinearEq}
  (i\partial_t + h(D))u = F(t,x), \qquad u(0) = f(x),
\end{equation}
which has the representation formula
\begin{equation}\label{RepForm}
  u(t) = \mathcal{W}(t) f + \int_0^t \mathcal{W}(t-t')F(t') \, dt',
\end{equation}
 where  $\mathcal{W}(t) = e^{ith(D)}$ is the associated solution
group. 
We define the associated $X^{s,b}_{p}$ spaces with the norm
\begin{align*}
   \norm{u}_{X^{s,b}_{p} }&= \bignorm{\angles{\xi}^s
\angles{-\tau+h(\xi)}^b \widetilde
u(\tau,\xi)}_{L^{p'}_{\tau\xi}}.
          \end{align*}
 In the special cases, $h(\xi)=\pm\abs{\xi}$, we use the notation $ X^{s,b}_{p,\pm}$. We also define
 $ \widehat{L^p_{tx} }=X^{0,0}_{p} $.    
The restriction of these spaces to a time slab $S_T=[-T, T]\times
\R$ are denoted $X_p^{s,b}(S_T)$,  defined in the usual way by the norm 
$$
  \norm{u}_{X^{s,b}_{p} (S_T)}=\inf \left\{ \norm{u'}_{X^{s,b}_{p} }:  u' =u \ \text{on} \ S_T  \right\}.
$$

For $b>\frac1p$, we have the embedding
\begin{equation}\label{X-H}
X_{p}^{s,b} \subset  C\left(\R; \widehat{H_p^s} \right)
  \end{equation}
where $C$ depends only on $b$. Indeed,
\begin{align*}
 \norm{u(t)}^{p'}_{\widehat{H_p^s} }=  \int \angles{\xi}^{p's} \abs{ \widehat{u} (t, \xi)}^{p'} \, d\xi=\int \angles{\xi}^{p's} \Bigabs{ \int e^{it\tau} \widetilde{u} (\tau, \xi) d\tau}^{p'} \, d\xi,
\end{align*}
but by H\"{o}lder inequality
\begin{align*}
 \Bigabs{ \int e^{it\tau} \widetilde{u} (\tau, \xi) d\tau}\le  \left( \int  \angles{-\tau+h(\xi)}^{-pb} \, d\tau \right)^{\frac1p} \left( \int  \angles{-\tau+h(\xi)}^{p'b}  \abs{ \widetilde{u} (\tau, \xi) }^{p'} \, d\tau \right)^{\frac1p'},
\end{align*}
where the first integral on the right is bounded since $bp>1$. A combination of these estimates will imply \eqref{X-H}. 

We also need the following estimate for the solution \eqref{RepForm} of the  inhomogeneous linear equation \eqref{LinearEq}. The proof is included in the last section by modifying  the proof of Lemma 5 in ~\cite{dfs2007}.
\begin{lemma}
\label{LemmaLinearEst}
Let $1<p<\infty$, $\frac1p<b\le 1$, $0\le \varepsilon \le 1-b$ and $0<T<1$. 
Assume $f \in \widehat{H_p^s} $ and $F \in X_p^{s,
b-1+\varepsilon}(S_T)$. Then $u$  in \eqref{RepForm} satisfies
$$
\norm{u}_{ X_{p}^{s,b}(S_T)  }\le C \norm{f}_{\widehat{H_p^s} }+ C T^\varepsilon
\norm{F}_{X_{p}^{s,b-1+\varepsilon} (S_T)},
$$
where $C$ depends on $b$.
\end{lemma}

\section{Reduction of Theorem \ref{Thm} to null form estimates}
Set $\frac1p=1-2\varepsilon $ for $0<\varepsilon \le \frac14$. By persistence of higher regularity argument it suffices to prove Theorem \ref{Thm} for $s=\frac1p$.

We shall then iterate the solutions to \eqref{MonopEq2-v3}--\eqref{Data2}  in 
$$
u_\pm, v_\pm \in X_{p, \pm }^{\frac1p, \frac1p+\varepsilon}(S_T).
$$
By a standard argument, using Lemma \ref{LemmaLinearEst}, the local existence problem of Theorem \ref{Thm} reduces to proving the bilinear estimate (the estimate for $N(v,u)$ is symmetrical)
\begin{equation}
\label{ReducBiEstimate}
\norm{ \mathcal{P}_\pm(D) N(u,v)}_{ X_{p, \pm }^{\frac1p, 0} } \lesssim \left( \norm{u_+}_{X_{p,+}^{\frac1p, \frac1p+\varepsilon}}+ \norm{u_-}_{X_{p, -}^{\frac1p, \frac1p+\varepsilon}}+\norm{v_+}_{X_{p, +}^{\frac1p, \frac1p+\varepsilon}} +\norm{v_-}_{X_{p, -}^{\frac1p, \frac1p+\varepsilon}}\right)^2.
\end{equation}
 
Writing $u=\mathcal{P}_+(D)u_++\mathcal{P}_-(D)u_-$ and $v=\mathcal{P}_+(D)v_++\mathcal{P}_- (D)v_-$,  we get 
$$
 N(u,v)=\sum_{\pm_1, \pm_2}\left(\begin{matrix}
 \frac12\left(\mathcal{P}_{\pm_1}(D)u_{\pm_1}\cdot \mathcal{P}_{\pm_2}(D)v_{\pm_2}-\mathcal{P}_{\pm_2}(D)v_{\pm_2}\cdot
  \mathcal{P}_{\pm_1}(D)u_{\pm_1}\right)\\
 \mathcal{P}_{\mp_1}(D) \beta u_{\pm_1}\cdot\mathcal{P}_{\pm_2}(D) u_{\pm_2}
\end{matrix}\right ),
 $$
 where in the second row we used the identity $ \beta \mathcal{P}_\pm(D)=  \mathcal{P}_\mp(D)\beta$; the signs ${\pm_1}$ and $ {\pm_2}$  are independent.
 Now, noting the property 
  $$\norm{u(-t, x)}_{X_{p, \pm}^{s, b}}= \norm{u(t, x)}_{X_{p, \mp}^{s, b}}, \quad \norm{u(t, -x)}_{X_{p, \pm}^{s, b}}= \norm{u(t, x)}_{X_{p, \pm}^{s, b}},$$
  and that $\mathcal{P}_\pm(\xi)$ is bounded,  \eqref{ReducBiEstimate} will reduce to proving
  \begin{equation}
\label{ReducBiEstimate-2}
\norm{\mathcal{P}_+(D) w\cdot \mathcal{P}_\pm(D) z }_{ X_{p}^{\frac1p, 0} } \lesssim \norm{w}_{X_{p,+}^{\frac1p, \frac1p+\varepsilon}}\norm{z}_{X_{p,\mp}^{\frac1p, \frac1p+\varepsilon}},
\end{equation}
for $w,z \in \mathcal{S}(\R^{1+2})$ taking values in $\mathfrak{g}\times \mathfrak{g}$.

 The key observation in the proof of \eqref{ReducBiEstimate-2} is the bilinear terms
   $\mathcal{P}_+(D) w\cdot \mathcal{P}_\pm(D) z$ are  null forms as shown in 
   ~\cite{cb2012}. This can be seen by taking their space-time Fourier transform. Indeed,
\begin{equation}\label{FT}
\widetilde{[\mathcal{P}_+(D) w\cdot \mathcal{P}_\pm(D) z]}(\tau,\xi)
  =\int_{\R^{1+2}} \mathcal{P}_{\pm}(\xi-\eta) \mathcal{P}_{+}(\eta) \widetilde  w(\lambda, \eta) \cdot \widetilde z(\tau-\lambda,\xi-\eta) \, d\lambda \, d\eta,
\end{equation}
 where the symbol  satisfies the estimate (see Lemma 2.2. in ~\cite{cb2012})
 \begin{equation}
 \label{NullSymbolEst}
\abs{ \mathcal{P}_{\pm}(\xi-\eta) \mathcal{P}_{+}(\eta) }\lesssim \theta(\eta, \mp( \xi-\eta)).
 \end{equation}
 The angles on the right hand side quantifies the null structure (see Lemma 2 and Remark 2  in ~\cite{dfs2007} ).

 So in view of \eqref{FT} and \eqref{NullSymbolEst}, the proof of \eqref{ReducBiEstimate-2} essentially reduces to
 \begin{equation}
\label{KeyBiEstimate}
\norm{ \mathcal{Q}_{\pm}   (\phi, \psi) }_{ X_{p }^{\frac1p, 0}} \lesssim \norm{\phi}_{X_{p, +}^{\frac1p, \frac1p+\varepsilon}}\norm{\psi}_{X_{p, \pm}^{\frac1p, \frac1p+\varepsilon}},
\end{equation}
where 
\begin{equation}\label{NullFormTerm}
\widetilde{ \mathcal{Q}_{ \pm  }(\phi, \psi)}(\tau,\xi)
  =\int_{\R^{1+2}} \theta(\eta,\pm(\xi-\eta)) \widetilde
\phi(\lambda,\eta)
  \widetilde \psi (\tau-\lambda,\xi-\eta) \, d\lambda \, d\eta.
\end{equation}
for $\phi, \psi:  \R^{1+2}\rightarrow \R$ such that $\widetilde
\phi, \ \widetilde\psi >0$.

So everything boils down to proving \eqref{KeyBiEstimate}. We identify two cases  where the product $\phi \psi $ has Fourier support contained in either of the following sets:
\begin{enumerate}[(I)]
\item $\{\xi: \abs{\xi}< 1 \}$ , low frequency  case
\item $\{\xi: \abs{\xi}\ge 1 \}$,  high frequency case.
\end{enumerate}
We give the proof of \eqref{KeyBiEstimate} in both of these cases in the following two Sections.
\section{Proof of \eqref{KeyBiEstimate} in the case of (I)}
Let $\chi\in C^\infty_c(\R^2)$ such that $\chi(\xi)=1$ on the set $\{\xi: \abs{\xi}< 1 \}.$  Then  we estimate
\begin{align*}
\text{l.h.s of  \eqref{KeyBiEstimate} }\le \norm{\widetilde{\phi\psi }}_{L^{p'}_{\tau \xi}}=
\norm{\chi\widetilde{\phi \psi }}_{L^{p'}_{\tau \xi}}
\le \norm{\chi}_{L^{\frac1\varepsilon}_{\xi}}
 \norm{  \norm{\widetilde{\phi \psi } }_{L^{p'}_{\tau}  }}_{ L^{\frac1\varepsilon}_{ \xi} }\lesssim \norm{  \norm{\widetilde{\phi \psi } }_{L^{p'}_{\tau}  }}_{ L^{\frac1\varepsilon}_{ \xi} },
\end{align*}
where we used H\"{o}lder inequality in $\xi$ and the assumption that $\frac1{p'}=2\varepsilon$ (since  $\frac1{p}=1-2\varepsilon$).
Applying Young's inequality to the convolution $\widetilde{\phi \psi }=\widetilde{\phi } \underset{\tau,\xi}* \widetilde{\psi } $, first in $\tau$ and then in $\xi$, we obtain
\begin{align*}
 \norm{  \norm{\widetilde{\phi \psi } }_{L^{p'}_{\tau}  }}_{ L^{\frac1\varepsilon}_{ \xi} }\le   \norm{  \norm{\widetilde{\phi  }}_{L^{1}_{\tau}  } \underset{\xi}*   \norm{\widetilde{\psi }}_{ L^{p'}_{\tau}  }}_{ L^{\frac1\varepsilon}_{ \xi} }
 \le\norm{  \norm{\widetilde{\phi }}_{L^{1}_{\tau}  }  }_{ L^{\frac2{1+\varepsilon}}_{ \xi} }\norm{  \norm{\widetilde{\psi }}_{L^{p'}_{\tau}  }  }_{ L^{\frac2{1+\varepsilon}}_{ \xi} }.
\end{align*}

Now, by H\"{o}lder
\begin{align*}
 \norm{  \norm{\widetilde{\phi }}_{L^{1}_{\tau}  }  }_{ L^{\frac2{1+\varepsilon}}_{ \xi} } &\le   \norm{ \angles{\xi}^\frac1p \norm{ \angles{-\tau\pm \abs{\xi}}^{\frac1p+\varepsilon}\widetilde{\phi }}_{L^{p'}_{\tau}   }  \norm{ \angles{-\tau\pm \abs{\xi}}^{-(\frac1p+\varepsilon)} }_{L^{p}_{\tau}   }  }_{ L^{p'}_{ \xi} } \norm{ \angles{\xi}^\frac1p  }_{ L^{\frac2{1-3\varepsilon}}_{ \xi} } 
 \\
 &\lesssim \norm{ \angles{\xi}^\frac1p  \angles{-\tau\pm \abs{\xi}}^{\frac1p+\varepsilon}  \widetilde{\phi }}_{ L^{p'}_{ \tau\xi} } =\norm{\phi }_{X_{p, \pm}^{\frac1p, \frac1p+\varepsilon}},
\end{align*}
where we used the fact that 
$$ \norm{ \angles{-\tau\pm \abs{\xi}}^{-(\frac1p+\varepsilon)} }_{L^{p}_{\tau}   }  \lesssim 1, \quad \norm{ \angles{\xi}^\frac1p  }_{ L^{\frac2{1-3\varepsilon}}_{ \xi} } \lesssim 1, $$ since $(\frac1p+\varepsilon)p>1$ and $\frac1p\frac2{1-3\varepsilon}=2\frac{1-2\varepsilon}{1-3\varepsilon}>2$.

Similarly, 
\begin{align*}
 \norm{  \norm{\widetilde{\psi}}_{L^{p'}_{\tau}  }  }_{ L^{\frac2{1+\varepsilon}}_{ \xi} } &\le   \norm{ \angles{\xi}^\frac1p \norm{\widetilde{ \psi }}_{L^{p'}_{\tau}   }  }_{ L^{p'}_{ \xi} } \norm{ \angles{\xi}^\frac1p  }_{ L^{\frac2{1-3\varepsilon}}_{ \xi} } 
 \\
 &\lesssim \norm{ \angles{\xi}^\frac1p \widetilde{\psi }}_{ L^{p'}_{ \tau\xi} } =\norm{ \psi }_{X_{p, \pm}^{\frac1p, 0}}.
\end{align*}

\section{Proof of \eqref{KeyBiEstimate}  in the case of (II)}
The angles in \eqref{NullFormTerm} satisfy the following estimates (see eg. ~\cite{dfs2006}, ~\cite{cb2012}):
\begin{equation}\label{AngleEst}
\theta(\eta, \xi-\eta)\approx \frac{r_+^\frac12}{\min(\abs{\eta}, \abs{\xi-\eta})^\frac12}, \qquad \theta(\eta, -(\xi-\eta))\approx \frac{\abs{\xi}^\frac12r_-^\frac12}{\abs{\eta}^\frac12\abs{\xi-\eta}^\frac12},
\end{equation}
where
\begin{align*}
r_+= \abs{\eta}+\abs{\xi-\eta}-\abs{\xi}, \qquad r_-=\abs{\xi}-\abs{\abs{\eta}- \abs{\xi-\eta}}.
\end{align*}

So in view of \eqref{NullFormTerm} and \eqref{AngleEst}, the estimate \eqref{KeyBiEstimate} reduces to 
\begin{align}
\label{BiEst++1}
\norm{\abs{D}^\frac1p \mathcal{R}_+^\frac12(\phi, \psi )}_{ \widehat{L^p_{tx} }}&\lesssim   \norm{\abs{D}^\frac1p \phi }_{X^{0, \frac1p+\varepsilon}_{p,+} }\norm{\abs{D}^{\frac1p+\frac12} \psi }_{X^{0, \frac1p+\varepsilon}_{p,+} },\\
\label{BiEst+-1}
\norm{\abs{D}^{\frac1p+\frac12} \mathcal{R}_-^\frac12(\phi , \psi )}_{ \widehat{L^p_{tx} }}&\lesssim   \norm{\abs{D}^{\frac1p+\frac12} \phi }_{X^{0,\frac1p+\varepsilon}_{p,+} }\norm{\abs{D}^{\frac1p+\frac12} \psi}_{X^{0,\frac1p+\varepsilon}_{p,-} },
\end{align}
where
\begin{equation*}
\widetilde{ \mathcal{R}_\pm^\frac12(\phi , \psi )}(\tau,\xi)
  =\int_{\R^{1+2}} r_\pm^\frac12 \ \widetilde
\phi (\lambda,\eta)
  \widetilde \psi  (\tau-\lambda,\xi-\eta) \, d\lambda \, d\eta.
\end{equation*}

By transfer principle (see eg. Proposition A.2 in ~\cite{ga2013}), the estimates \eqref{BiEst++1} and \eqref{BiEst+-1} further reduce to 
\begin{align}
\label{BiEst++2}
\norm{\abs{D}^\frac1p \mathcal{R}_+^\frac12\left(e^{it\abs{D} }f, e^{it\abs{D} }g\right)}_{ \widehat{L^p_{tx} }}&\lesssim   \norm{\abs{D}^\frac1p f}_{\widehat{L^p_{x} }}\norm{\abs{D}^{\frac1p+\frac12} g}_{\widehat{L^p_{x} }},\\
\label{BiEst+-2}
\norm{\abs{D}^{\frac1p+\frac12} \mathcal{R}_-^\frac12\left(e^{it\abs{D} }f, e^{-it\abs{D} }g\right)}_{ \widehat{L^p_{tx} }}&\lesssim   \norm{\abs{D}^{\frac1p+\frac12} f}_{\widehat{L^p_{x} } }\norm{\abs{D}^{\frac1p+\frac12}                                              g}_{\widehat{L^p_{x} } }.
\end{align}
These reduce to estimating
\begin{equation}
\label{ReducBiEst}
\norm{\mathcal{B}_{+, \pm}(f, g)}_{ \widehat{L^p_{tx} }}\lesssim   \norm{ f}_{\widehat{L^p_{x} }}\norm{g}_{\widehat{L^p_{x} }},
\end{equation}
where
\begin{align*}
\widetilde{\mathcal{B}_{+, +}(f, g)}(\tau,\xi)
  =\int_{\R^{2}}  \frac{ \abs{\xi}^\frac1p (\tau-\abs{\xi})^\frac12 } 
    {\abs{\eta}^\frac1p \abs{\xi-\eta}^{\frac1p+\frac12} }  \hat{f}(\eta)
  \hat{g}(\xi-\eta) \delta(\tau-\abs{\eta}-\abs{\xi-\eta})\,  d\eta
  \\
  \widetilde{\mathcal{B}_{+, -}(f, g)}(\tau,\xi)
  =\int_{\R^{2}}  \frac{ \abs{\xi}^{\frac1p+\frac12} \abs{\abs{\tau}-\abs{\xi}}^\frac12 } 
    {\abs{\eta}^{\frac1p+\frac12} \abs{\xi-\eta}^{\frac1p+\frac12} }  \hat{f}(\eta)
  \hat{g}(\xi-\eta) \delta(\tau-\abs{\eta}+\abs{\xi-\eta})\,  d\eta.
\end{align*}

\subsection{Proof of \eqref{ReducBiEst} for $B_{+,+}$}
By H\"{o}lder inequality 
$$
\abs{ \widetilde{\mathcal{B}_{+, +}(f, g) }(\tau,\xi)  }\le \abs{ I(\tau, \xi) }^\frac1p 
\left\{   \int_{\R^{2}}  \abs{\hat{f}(\eta)}^{p'}\abs{ \hat{g}(\xi-\eta)}^{p'} \delta(\tau-\abs{\eta}-\abs{\xi-\eta}) \, d\eta \right\}^{\frac{1}{p'}}
$$
where 
$$
I(\tau, \xi)=  \abs{\xi} (\tau-\abs{\xi})^\frac p2 \int_{ \R^{2} }  \frac{ \delta(\tau-\abs{\eta}-\abs{\xi-\eta})  } 
    { \abs{\eta} \abs{\xi-\eta}^{1+\frac p2}  } \,  d\eta.
$$
By Proposition 4.3 in~\cite{fk2000}, 
$$
\int_{ \R^{2} }  \frac{ \delta(\tau-\abs{\eta}-\abs{\xi-\eta})  } 
    { \abs{\eta} \abs{\xi-\eta}^{1+\frac p2}  } \,  d\eta \approx \frac{1}{\tau (\tau-\abs{\xi})^\frac p2}\le  \frac{1}{\abs{\xi} (\tau-\abs{\xi})^\frac p2},
$$
where  we used the fact that $p>1$ (in the first inequality) and $\tau>\abs{\xi}$ (in the second inequality). Thus, we get
$$
I(\tau, \xi)\lesssim 1 \quad \text{for all} \ \tau, \xi,
$$
which in turn implies
\begin{align*}
\int_{ \R^{1+2} }  \abs{ \widetilde{\mathcal{B}_{+, +}(f, g) }(\tau,\xi)  }^{p'}\, d\tau \, d\xi &\lesssim \int_{\R^{1+5}}  \abs{\hat{f}(\eta)}^{p'}\abs{ \hat{g}(\xi-\eta)}^{p'} \delta(\tau-\abs{\eta}-\abs{\xi-\eta}) \, d\tau \, d\xi\, d\eta
\\
&= \int_{\R^{2}}  \abs{\hat{f}(\eta)}^{p'} d\eta  \int_{\R^{2}}  \abs{\hat{g}(\eta)}^{p'} d\eta,
\end{align*}
where to get the equality we used the fact that $\int_{\R}  \delta(\tau-\abs{\eta}-\abs{\xi-\eta}) \, d\tau =1$. This proves \eqref{ReducBiEst} in the $B_{+,+}$ case.

\subsection{Estimate for $B_{+,-}$}
Again, by H\"{o}lder inequality 
$$
\abs{ \widetilde{\mathcal{B}_{+, -}(f, g) }(\tau,\xi)  }\le \abs{ J(\tau, \xi) }^\frac1p 
\left\{   \int_{\R^{2}}  \abs{\hat{f}(\eta)}^{p'}\abs{ \hat{g}(\xi-\eta)}^{p'} \delta(\tau-\abs{\eta}+\abs{\xi-\eta}) \, d\eta \right\}^{\frac{1}{p'}}
$$
where 
$$
J(\tau, \xi)=  \abs{\xi}^{1+\frac p2}  \Bigabs{ \abs{\xi}-\abs{\tau} } ^\frac p2 \int_{ \R^{2} }  \frac{ \delta(\tau-\abs{\eta}+\abs{\xi-\eta})  } 
    { \abs{\eta}^{1+\frac p2} \abs{\xi-\eta}^{1+\frac p2}  } \,  d\eta.
$$
By a similar argument as in the preceding subsection, it then suffices to show 
$$
J(\tau, \xi)\lesssim 1 \quad \text{for all} \ \tau, \xi.
$$
To do so, we split the integral
$$
 \int_{ \R^{2} }  \frac{ \delta(\tau-\abs{\eta}+\abs{\xi-\eta})  } 
    { \abs{\eta}^{1+\frac p2} \abs{\xi-\eta}^{1+\frac p2}  } \,  d\eta=J_1+ J_2,
$$
where 
\begin{align*}
J_1&=\int_{\abs{\eta}+\abs{\xi-\eta} \le 2\abs{\xi}  }  \frac{ \delta(\tau-\abs{\eta}+\abs{\xi-\eta})  } 
    { \abs{\eta}^{1+\frac p2} \abs{\xi-\eta}^{1+\frac p2}  }\, d\eta,
    \\
    J_2&=\int_{\abs{\eta}+\abs{\xi-\eta} \ge 2\abs{\xi}  }  \frac{ \delta(\tau-\abs{\eta}+\abs{\xi-\eta})  } 
    { \abs{\eta}^{1+\frac p2} \abs{\xi-\eta}^{1+\frac p2}  }\, d\eta.
\end{align*}
For $J_1$, we apply Lemma 4.5 in ~\cite{fk2000} (using the fact that $p>1$) to get the estimate
$$
J_1\approx \frac{1}{\abs{\xi}^{1+\frac p2}  \Bigabs{ \abs{\xi}-\abs{\tau} } ^\frac p2 },
$$
which in turn implies $J(\tau, \xi)\lesssim 1 \ \text{for all} \ \tau, \xi.
$.

Next, we consider $J_2$.  Applying Lemma 4.4 in ~\cite{fk2000} with
 $$F(\abs{\eta}, \abs{\xi-\eta})=  \abs{\eta}^{-1-\frac p2} \abs{\xi-\eta}^{-1-\frac p2} \mathbbm{1}_{\{  \abs{\eta}+\abs{\xi-\eta} \ge 2\abs{\xi}\}}.$$ 
 we get 
 $$
 J_2\approx \Bigabs{ \abs{\xi}^2-\tau^2}^{-\frac12} \int_1^\infty \Bigabs{\frac{ x\abs{\xi}+\tau}2}^{-1-\frac p2}\Bigabs{\frac{ x\abs{\xi}-\tau}2}^{-1-\frac p2}\mathbbm{1}_{\{ x\ge 2\}}
\Bigabs{\abs{\xi}^2x^2-\tau^2}(x^2-1)^{-\frac12} \, dx .
$$
Thus
 \begin{align*}
 J_2&\approx \Bigabs{\abs{\xi}^2-\tau^2}^{-\frac12} \abs{\xi}^{-p}\int_2^\infty \Bigabs{x^2-\frac{\tau^2}{\abs{\xi}^2}}^{-\frac p2}(x^2-1)^{-\frac12}\, dx ,
 \\
 &\lesssim  \Bigabs{\abs{\xi}^2-\tau^2}^{-\frac12} \abs{\xi}^{-p}\int_2^\infty x^{-p-1}\, dx 
 \\
 &\lesssim \abs{\xi}^{-p-\frac12}\Bigabs{\abs{\xi}-\abs{\tau}}^{-\frac12} ,
 \end{align*}
 where we used the fact that $\tau\le \abs{\xi}$  and  $\int_2^\infty x^{-p-1}\, dx \lesssim 1$.
Using the estimate for $J_2$ in $J$, we get 
$$J(\tau, \xi)\lesssim \left\{  \frac{  \Bigabs{ \abs{\xi}-\abs{\tau} } }{ \abs{\xi} } \right\}^{\frac {p-1}2}\le 1  \quad \text{for all} \ \tau, \xi. $$

\section{Proof of Lemma \ref{LemmaLinearEst}}
Let  $\rho \in C_c^\infty([-1,1])$ be a cut-off function such that
$\rho(t) = 1$ for $\abs{t} \le 1$ and $\rho(t) = 0$ for $\abs{t} \ge
2$. Set $\rho_T(t) = \rho(t/T)$. 

  The homogeneous part of the solution satisfies the estimate 
$$
\norm{\mathcal{W}(t)f}_{ X_{p}^{s,b}(S_T)  }\le 
\norm{\rho_T(t) \mathcal{W}(t)f}_{ X_{p}^{s,b}}=\norm{\rho_T}_{\widehat{H_{p}^{b}}}\norm{f}_{\widehat{H_{p}^{s}}}\le C_\rho \norm{f}_{\widehat{H_{p}^{s}}}.
$$

Next, consider the inhomogeneous part of the solution. Write $$v(t)= \int_0^t \mathcal{W}(t-t')
F(t')dt'.$$ Extend $F$ by zero outside $S_T$.
 Taking
Fourier transform in space,
\begin{equation*}\label{vFourier}
  \widehat v(t,\xi)
  =
  \int_0^t e^{ i(t-t')h(\xi)} \widehat F(t',\xi) \, dt'
  \simeq
  \int \frac{e^{it\lambda}-e^{ith(\xi)}}{i(\lambda-h(\xi))} \widetilde F(\lambda,\xi) \, d\lambda
\end{equation*}
and then also in time,
\begin{align*}
  \widetilde v(\tau,\xi)
  &=
  \int \frac{\delta(\tau-\lambda)-\delta(\tau- h(\xi))}{i(\lambda- h(\xi))} \widetilde F(\lambda,\xi) \,
  d\lambda\\
 & =
  \frac{\widetilde F(\tau,\xi)}{i(\tau- h(\xi))}
  - \delta(\tau- h(\xi)) \int \frac{\widetilde F (\lambda,\xi)}{i(\lambda - h(\xi))}  \, d\lambda.
\end{align*}

Now split $F = F_1 + F_2$ corresponding to the Fourier domains
$ \left\{T\abs{\lambda -h(\xi)} \lesssim 1\right\}$ and $ \left\{T\abs{\lambda -h(\xi)} \gg 1\right\}$ respectively. Write $v=v_1+v_2$ accordingly. Expand
$$
  \widehat{v_1}(t,\xi)
  =
  e^{ ith(\xi)} \sum_{n=1}^\infty \int \frac{\left[it(\lambda - h(\xi))\right]^n}{n! i(\lambda - h(\xi))} \mathbbm{1}_{ \left\{T\abs{\lambda -h(\xi)} \lesssim 1\right\}} \widetilde F(\lambda,\xi) \, d\lambda
$$
hence
\begin{equation*}\label{v1}
  v_1(t) =
  \sum_{n=1}^\infty \frac{t^n}{n!} \mathcal{W}(t) f_n
\end{equation*}
where
\begin{equation*}\label{fn}
  \widehat{f_n}(\xi) = \int \left[i(\lambda -h(\xi))\right]^{n-1}
  \mathbbm{1}_{ \left\{T\abs{\lambda -h(\xi)} \lesssim 1\right\}} \widetilde F(\lambda,\xi) \, d\lambda.
\end{equation*}
Then
\begin{align*}
  \norm{v_1}_{X_p^{s,b}(S_T)}
  \le
  \norm{\rho_T v_1}_{X_p^{s,b}}
 & \le
  \sum_{n=1}^\infty \frac{T^n}{n!} \norm{\left(\frac{t}{ T}\right)^n \rho_T(t)W(t)f_n}_{X_p^{s,b}}
  \\
  &\le
  \sum_{n=1}^\infty \frac{T^n}{n!} \norm{\left(\frac{t}{ T}\right)^n \rho_T(t)}_{\widehat{H_{p}^{b}}}
  \norm{f_n}_{\widehat{H_{p}^{s}}}.
\end{align*}

Using H\"{o}der inequality with respect to the variable $\lambda$, we get 
\begin{equation}\label{Estfn}
 \norm{f_n}_{ \widehat{H_{p}^{s}} }\lesssim  T^{b+\varepsilon-n-\frac1p}
  \norm{F}_{X_p^{s,b-1+\varepsilon} } . 
\end{equation}
On the other hand, we can estimate
\begin{equation}
\label{EstRho}
 \norm{\left(\frac{t}{ T}\right)^n \rho_T}_{\widehat{H_{p}^{b}}}\le C  T^{\frac1p-b} \norm{t^n \rho}_{\widehat{H_{p}^{b}}} \le C  T^{\frac1p-b} n2^{n}.
\end{equation}
Indeed, 
\begin{align*}
\norm{\left(\frac{t}{ T}\right)^n \rho_T}_{\widehat{H_{p}^{b}}}^{p'}&=T^{-np'}\int \angles{\tau}^{p' b}  \abs{ \widehat{\rho_T}^{(n)}(\tau) }^{p'}\, d\tau
\\
&=T^{p'} \int \angles{\tau}^{p' b}  \abs{ \widehat{\rho}^{(n)}(T\tau) }^{p'}\, d\tau
\\
&\le CT^{p'-p'b-1} \int \angles{\tau}^{p' b}  \abs{ \widehat{\rho}^{(n)}(\tau) }^{p'}\, d\tau
\\
&=CT^{p'-p'b-1} \norm{t^n \rho}_{\widehat{H_{p}^{b}}}^{p'},
\end{align*}
from which we get the first inequality in \eqref{EstRho} by taking the $p'$-th root. Whereas the second inequality can be estimated using the support assumption of $\rho$ as 
$$
 \norm{t^n \rho}_{\widehat{H_{p}^{b}}}\le C \norm{t^n \rho}_{H^1}\le C n 2^n \norm{\rho}_{H^1}.
$$

So in view of \eqref{Estfn} and \eqref{EstRho}, we have
\begin{align*}
  \norm{v_1}_{X_p^{s,b}(S_T)}
    &\le CT^{\varepsilon}
  \left(\sum_{n=1}^\infty \frac{n2^{n}}{n!}\right)
  \norm{F}_{X^{s,b-1+\delta}}\\
  &\le C T^{\varepsilon}
    \norm{F}_{X^{s,b-1+\varepsilon}}.
\end{align*}

It remains to prove the estimate for $v_2$. We split $v_2 = w_1-w_2$, where
\begin{align*}
   \widetilde w_1(\tau,\xi)
 & =
  \frac{\mathbbm{1}_{  \left\{T\abs{\lambda -h(\xi)} \gg 1\right\} }\widetilde F(\tau,\xi)}{i(\tau-h(\xi))},
  \\
   \widetilde w_2(\tau,\xi)
&  =
  \delta(\tau-h(\xi)) \widehat g(\xi),
\end{align*}
for 
$$
  \widehat g(\xi)
  =
  \int \frac{\mathbbm{1}_{  \left\{T\abs{\lambda -h(\xi)} \gg 1\right\}} \widetilde F(\lambda,\xi) }{i(\lambda-h(\xi))}  \, d\lambda.
$$
Obviously,
$$
  \norm{w_1}_{X_p^{s,b}(S_T)}\lesssim T^\varepsilon  \norm{F}_{X_p^{s,b-1+\varepsilon}}.
$$
On the other hand,
$$
\norm{w_2}_{X_p^{s,b}(S_T)}\le\norm{\rho_T w_1}_{X^{s,b}}\le
\norm{\rho_T}_{\widehat{H_{p}^{b}}}\norm{g}_{\widehat{H_{p}^{s}}}.
$$
Now, one can easily show 
$$
\norm{\rho_T}_{\widehat{H_{p}^{b}}}\le C T^{\frac1p-b} \norm{\rho}_{\widehat{H_{p}^{b}}}\le C_\rho T^{\frac1p-b} 
$$
and using H\"{o}lder inequality, we obtain
$$
\norm{g}_{\widehat{H_{p}^{s}}}\le T^{b+\varepsilon-\frac1p}\norm{F}_{X^{s,b-1+\varepsilon}}.
$$
A combination of these estimates gives the desired estimate for $w_2$.

\end{document}